\documentclass[11pt, reqno]{amsart}
\usepackage{latexsym}
\usepackage{amssymb}
\usepackage{amsthm}
\usepackage{amsmath}
\usepackage{graphicx}
\usepackage{color}
\allowdisplaybreaks   

\newtheorem{theorem}{Theorem}[section]

\newtheorem{corollary}[theorem]{Corollary}

\newtheorem{lemma}[theorem]{Lemma}

\newtheorem{definition}[theorem]{Definition}

\newtheorem{proposition}[theorem]{Proposition}

\newtheorem*{Remark}{Remark}
\numberwithin{equation}{section} 
\makeatletter
\let\c@equation\c@theorem 
\makeatother

\DeclareMathOperator{\Hom}{Hom}

\DeclareMathOperator{\Ext}{Ext}

\DeclareMathOperator{\HH}{HH}
\DeclareMathOperator{\Ho}{H}
\newcommand{\bk}{\mathbf{k}}
\newcommand{\kG}{\bk G}

\DeclareMathOperator{\res}{res}
\DeclareMathOperator{\cor}{cor}

\begin{document}
\begin{abstract}

\small{We show that the Tate-Hochschild cohomology ring $\widehat{\HH}^*(RG,RG)$ of a finite group algebra $RG$ is isomorphic to a direct sum of the Tate cohomology rings of the centralizers of conjugacy class representatives of $G$. Moreover, our main result provides an explicit formula for the cup product in $\widehat{\HH}^*(RG,RG)$ with respect to this decomposition. As an example, this formula helps us to compute the Tate-Hochschild cohomology ring of the symmetric group $S_3$ with coefficients in a field of characteristic $3$. 
\vspace{-0.3em}
}

\end{abstract}

\title[The Tate-Hochschild Cohomology Ring of a Group Algebra]
{The Tate-Hochschild \\
Cohomology Ring of a Group Algebra}

\author{Van C. Nguyen}
\address{Department of Mathematics\\Texas A\&M University \\College Station, TX 77843}
\email{vcnguyen@math.tamu.edu}

\maketitle

\section{Introduction}
\label{sec:intro}

The theory of group cohomology is a well-studied yet ongoing research area. It has many applications to other areas such as representation theory, algebraic geometry, and commutative algebra. It is well-known that the Hochschild cohomology of a group algebra agrees with its usual cohomology, in the following sense. For an arbitrary commutative ring $R$ and a group $G$, if $M$ is an $RG$-bimodule, then it may be regarded as a left $RG$-module via the diagonal action $g \cdot m = gmg^{-1}$, where $g \in G, m \in M$. Conversely, any left $RG$-module $M$ may be considered as an $RG$-bimodule by letting $RG$ act trivially on the right. Together with the Eckmann-Shapiro Lemma, this shows that the Hochschild cohomology ring, $\HH^*(RG,M):= \bigoplus_{n \geq 0} \Ext^n_{RG \otimes_R RG^{op}}(RG,M)$ with coefficients in an $RG$-bimodule $M$, is the same as the usual group cohomology ring, $\Ho^*(G,M):= \bigoplus_{n \geq 0} \Ext^n_{RG}(R,M)$ with coefficients in $M$ under the diagonal action. In particular, by considering $RG$ as its own bimodule, $\HH^*(RG,RG)$ is isomorphic to $\Ho^*(G,RG)$, where $RG$ is a left $RG$-module via conjugation. From this identification and the Eckmann-Shapiro Lemma, one can prove that $\HH^*(RG,RG)$ may be decomposed as a direct sum of the cohomology of the centralizers of conjugacy class representatives of $G$ (\cite{Ben2}, Theorem~2.11.2). In 1999, Siegel and Witherspoon then described a formula for the products in $\HH^*(RG,RG)$ in terms of this additive decomposition \cite{SW}. When $G$ is abelian, Cibils and Solotar proved that the Hochschild cohomology ring of $G$ is (isomorphic to) the tensor product over $R$ of $RG$ and its usual cohomology ring \cite{CiSo}. 

In 1952, John Tate introduced a group cohomology theory that is based on complete resolutions, and hence, expands the study of group cohomology to negative degrees \cite{Tate}. Tate cohomology exploits the similarities between the usual homology and cohomology of a finite group $G$. A lot of study about this new cohomology has been done and can be found in (\cite{Brown}, Ch.~VI) or (\cite{CaEi}, Ch.~XII). In particular, in 1992, Benson and Carlson proved that for a finite group $G$ and a field $\bk$ of characteristic $p >0$, if the depth of $\Ho^*(G,\bk)$ is greater than one, products of elements in negative cohomology are zero \cite{BeCa}. More recently, Bergh and Jorgensen combined the notions of Hochschild cohomology and Tate cohomology to form the Tate-Hochschild cohomology \cite{BeJo}. In this paper, we explore the structure of the Tate-Hochschild cohomology of a group algebra and generalize some known results about group cohomology to its negative degrees. 

The material is organized as follows. In Section~\ref{sec:Tate RG}, we recall the definition and properties of Tate cohomology. Let $M$ be a left $RG$-module. Then for all $n \in \mathbb{Z}$:
\[\widehat{\Ho}^n(G,M) := \widehat{\Ext}^n_{RG}(R,M). \]
The readers who are familiar with Tate cohomology can skip this section. In Section~\ref{sec:Tate-Hochschild RG}, we specialize the definition of Tate-Hochschild cohomology in \cite{BeJo} to the group algebra. Let $M$ be an $RG$-bimodule. Then for all $n \in \mathbb{Z}$:
\[\widehat{\HH}^n(RG,M) := \widehat{\Ext}^n_{RG \otimes_R RG^{op}}(RG,M). \]
In this section, we also briefly present a cup product that makes the Tate-Hochschild cohomology
\[\widehat{\HH}^*(RG,RG) = \bigoplus_{n \in \mathbb{Z}} \widehat{\HH}^n(RG,RG)\]
become a graded ring \cite{Nguyen}. In Section~\ref{sec:relationship RG}, we reduce the Tate-Hochschild cohomology of $RG$ to its Tate cohomology using a similar argument as in the usual cohomology. In particular, we show that $\widehat{\HH}^*(RG,RG)$ is isomorphic to $\widehat{\Ho}^*(G,RG)$, where $RG$ is a left $RG$-module via conjugation. In addition, we state known Tate cohomology relations with subgroups, which will be useful in proving our main result. In Section~\ref{sec:decomp}, we introduce a generalization of the Tate-Hochschild cohomology ring of $RG$ by letting another finite group $H$ act on $G$ and consider $\widehat{\Ho}^*(H,RG)$. Let $g_1, \ldots, g_t \in G$ be representatives of the orbits of the action of $H$ on $G$. Let $H_i:=\text{Stab}_H(g_i)=\{h \in H \:|\: ^hg_i=g_i\}$ be the stabilizer of $g_i$. We show that $\widehat{\Ho}^*(H,RG)$ decomposes as:
\[\widehat{\Ho}^*(H,RG) \cong \displaystyle \bigoplus_i \widehat{\Ho}^*(H_i,R). \]
In the main Theorem~\ref{THdecomp product}, we describe the multiplicative structure of $\widehat{\Ho}^*(H,RG)$ by giving it a product formula in terms of this additive decomposition. Our work, which uses mainly the idea from \cite{SW}, is a straightforward generalization of the usual group cohomology results. We observe that the Tate-Hochschild cohomology ring $\widehat{\Ho}^*(RG,RG) \cong \widehat{\Ho}^*(G,RG)$ is a special case of $\widehat{\Ho}^*(H,RG)$ by letting $H=G$ act on itself by conjugation. Working with this generalization, we prove that the Tate-Hochschild cohomology of $G$ decomposes as a direct sum of the Tate cohomology of the centralizers of conjugacy class representatives of $G$ with coefficients in $R$. Moreover, the product formula reduces the computation of products in $\widehat{\HH}^*(RG,RG)$ to products within the Tate cohomology rings of certain subgroups of $G$. When $G$ is abelian, Cibils and Solotar's result in \cite{CiSo} is also generalized to the isomorphism $\widehat{\HH}^*(RG,RG) \cong RG \otimes_R \widehat{\Ho}^*(G,R)$ in Proposition~\ref{abelian RG}. 
In Section~\ref{sec:S3}, we describe the Tate-Hochschild cohomology ring of the symmetric group on three elements $S_3$ over a field $\bk$ of characteristic $3$, utilizing the product formula in Theorem~\ref{THdecomp product}.

\section{Tate Cohomology of a Group Algebra}
\label{sec:Tate RG}

Throughout this paper, we let $G$ be a finite group and $R$ be the ring of integers $\mathbb{Z}$ or a field $\bk$ of characteristic $p>0$ such that $p$ divides the order of $G$. All rings and algebras are assumed to possess a unit; all modules are assumed to be left modules unless stated otherwise; and tensor products will be over $R$ unless stated otherwise. If $R=\bk$ is a field, then $\kG$ is a finite dimensional, symmetric, and self-injective algebra over $\bk$ (\cite{Ben1}, Prop.~3.1.2), and hence, projective $\kG$-modules are the same as injective $\kG$-modules. If $G$ is acting on a set $X$, then we denote the action $^gx=gxg^{-1}$, for all $g \in G$ and $x \in X$. 

Let $M$ and $N$ be left $RG$-modules. Then for any $g\in G$, $m\in M$, $n\in N$, and $f\in \Hom_R(M,N)$, we observe some basic facts:
\begin{itemize}
\item $M \otimes N$ is a left $RG$-module via $g \cdot (m \otimes n)= (g \cdot m) \otimes (g \cdot n)$,
\item $\Hom_R(M,N)$ is a left $RG$-module via $(g \cdot f)(m)= g(f(g^{-1} \cdot m))$,
\item We may regard left $RG$-modules as right $RG$-modules via $m \cdot g= g^{-1} \cdot m$.
\end{itemize}  

The Tate cohomology for $RG$ is defined in both positive and negative degrees using the following general resolution:

\begin{definition}
Let $\mathcal{R}$ be a two-sided Noetherian ring. A \textbf{complete resolution} of a finitely generated $\mathcal{R}$-module $M$ is an exact complex $\mathbb{P}=\{\{P_i\}_{i \in \mathbb{Z}}, d_i: P_i \rightarrow P_{i-1}\}$ of finitely generated projective $\mathcal{R}$-modules such that:
\begin{enumerate}
 \item The dual complex $\Hom_\mathcal{R}(\mathbb{P},\mathcal{R})$ is also exact
 \item There exists a projective resolution $\mathbb{Q} \stackrel{\varepsilon}{\rightarrow} M$ of $M$ and a chain map $\mathbb{P} \stackrel{\varphi}{\rightarrow} \mathbb{Q}$ where $\varphi_n$ is bijective for $n \geq 0$ and $0$ for $n<0$.
\end{enumerate}
\end{definition}

We explicitly construct an $RG$-complete resolution of $R$ as follows. Let
 \[\cdots \xrightarrow{d_3} P_2 \xrightarrow{d_2} P_1 \xrightarrow{d_1} P_0 \xrightarrow{\varepsilon} R \rightarrow 0 \]
be an $RG$-projective resolution of $R$, where each $P_i$ is a finitely generated projective $RG$-module. Apply $\Hom_R(-,R)$ to get a dual sequence:
 \[0 \rightarrow R \rightarrow \Hom_R(P_0,R) \rightarrow \Hom_R(P_1,R) \rightarrow \Hom_R(P_2,R) \rightarrow \cdots, \] 
which is an exact sequence of $RG$-modules, since $R \cong \Hom_R(R,R)$ and $\Hom_R(-,R)$ is an exact functor when $R=\bk$ is a field (the case $R=\mathbb{Z}$ is shown in \cite{CaEi}, Prop.~XII.3.3). Splicing these two sequences together, one forms a doubly infinite sequence:
 \[\cdots \xrightarrow{d_3} P_2 \xrightarrow{d_2} P_1 \xrightarrow{d_1} P_0 \rightarrow \Hom_R(P_0,R) \rightarrow \Hom_R(P_1,R) \rightarrow \cdots\]
We use the notation $P_{-(n+1)} := \Hom_R(P_n,R)$. By definition, the above sequence is an $RG$-complete resolution of $R$:
 \[\mathbb{P}: \hspace{1 cm} \cdots \xrightarrow{d_3} P_2 \xrightarrow{d_2} P_1 \xrightarrow{d_1} P_0 \rightarrow P_{-1} \rightarrow P_{-2} \rightarrow \cdots\] 
 
Let $M$ be any (left) $RG$-module. Apply $\Hom_{RG}(-,M)$ to $\mathbb{P}$ and take the homology of this new complex, we obtain the \textbf{Tate cohomology} for $RG$: 
 \[\widehat{\Ho}^n(G,M) := \widehat{\Ext}^n_{RG}(R,M) = \Ho^n(\Hom_{RG}(\mathbb{P},M)), \text{ for all } n \in \mathbb{Z}.\]

Observe that in our context, naturally, the Tate (co)homology is independent of the complete resolution of $R$ (\cite{AvMa}, Theorem~5.2 and Lemma~5.3). One can see this by applying a complete chain map between two complete resolutions of $R$ in both positive and negative degrees. This is the generalized Comparison Theorem on complete resolutions (\cite{Brown}, Prop.~VI.3.3). 

We note that the Tate cohomology of $RG$ obtains a multiplicative structure as described in (\cite{Brown}, Section~VI.5) and (\cite{CaEi}, Sections~XII.4 and XII.5). In particular, if $M$ is a ring on which $G$ acts by automorphisms: $g\cdot(m_1m_2)=(g\cdot m_1)(g\cdot m_2)$, then $\widehat{\Ho}^*(G,M)$ becomes an associative graded ring. It is graded-commutative, in the sense that, for $\alpha \in \widehat{\Ho}^i(G,M)$ and $\beta \in \widehat{\Ho}^j(G,M)$, $\alpha \beta=(-1)^{ij}\beta \alpha$ (\cite{CaEi}, Props. XII.5.2 and XII.5.3). When $M=R$, we denote $\widehat{\Ho}^*(G):=\widehat{\Ho}^*(G,R)$. 

Moreover, for a left $RG$-module $M$, we recall the following properties of Tate cohomology (\cite{Brown}, Section~VI.5):
 \begin{enumerate}
 \renewcommand{\labelenumi}{(\alph{enumi})}  
  \item For all $n>0$, $\widehat{\Ho}^n(G,M) \cong \Ho^n(G,M)$.
    
  This follows from the construction of a complete resolution of $R$. The positive-degree component of a complete resolution of $R$ arises from a projective resolution of $R$ which is also used to form the usual cohomology groups $\Ho^n(G,M)$.\\
  
  \item The group $\widehat{\Ho}^0(G,M)$ is a quotient of $\Ho^0(G,M)$.

  \item The group $\widehat{\Ho}^{-1}(G,M)$ is the dually defined submodule of $H_0(G,M)$.
  
  These follow from the construction of complete resolutions.\\
  
  \item For all $n < -1$, we have isomorphisms: $\widehat{\Ho}^n(G,M) \cong H_{-(n+1)}(G,M)$.
  
  For any left finitely generated $RG$-module $P_t$, the map $\tau: \Hom_R(P_t,R) \otimes_{RG} M \rightarrow \Hom_{RG}(P_t,M)$ is an $RG$-module isomorphism (\cite{CaEi}, Section~XII.3 (5)). It is defined explicitly by: $\tau(f \otimes_{RG} m)(p)= \sum_{g\in G} f(g^{-1}p)gm$, where $m\in M, p \in P_t$, and $f \in \Hom_R(P_t,R)$. With $t=n<-1$, $P_n := \Hom_R(P_{-(n+1)},R)$. Since $P_{-(n+1)}$ is finitely generated, $\Hom(P_n,R)=\Hom(\Hom(P_{-(n+1)},R),R) \cong P_{-(n+1)}$. Hence, by the duality isomorphism, we have $P_{-(n+1)} \otimes_{RG} M \cong \Hom_{RG}(P_n,M)$, which induces the desired isomorphism on homology.\\
   
   \item If $0 \rightarrow M \rightarrow M' \rightarrow M'' \rightarrow 0$ is a short exact sequence of (left) $RG$-modules, then there is a doubly infinite long exact sequence of Tate cohomology groups:
 \[\cdots \rightarrow \widehat{\Ho}^n(G,M) \rightarrow \widehat{\Ho}^n(G,M') \rightarrow \widehat{\Ho}^n(G,M'') \rightarrow \widehat{\Ho}^{n+1}(G,M) \rightarrow \cdots \] 

This is true since the spliced complex $\mathbb{P}$ consists of projective modules, so the short exact sequence of modules $0 \rightarrow M \rightarrow M' \rightarrow M'' \rightarrow 0$ gives rise to a short exact sequence of complexes: 
\[0 \rightarrow \Hom_{RG}(\mathbb{P},M) \rightarrow \Hom_{RG}(\mathbb{P},M') \rightarrow \Hom_{RG}(\mathbb{P},M'') \rightarrow 0 \]
whose cohomology long exact sequence is the desired long exact sequence (\cite{AvMa}, Prop.~5.4). \\

 \item If $(N_j)_{j\in J}$ is a finite family of (left) $RG$-modules and $(M_i)_{i\in I}$ is any family of $RG$-modules, then there are natural isomorphisms, for all $n \in \mathbb{Z}$:
  \[\displaystyle \widehat{\Ext}^n_{RG}(\bigoplus_{j \in J}{N_j}, M) \cong \prod_{j \in J} \widehat{\Ext}^n_{RG}(N_j,M) \]
  \[\displaystyle \widehat{\Ext}^n_{RG}(N, \prod_{i \in I} M_i) \cong \prod_{i \in I} \widehat{\Ext}^n_{RG}(N,M_i). \]

The idea of this proof is similar to that of the usual $\Ext^n_{RG}$, using the analogous relation for $\Hom$ (\cite{AvMa}, Prop.~5.7).
 \end{enumerate}
 
In 2011, Bergh and Jorgensen introduced the notion of Tate-Hochschild cohomology which is the extended Hochschild cohomology using complete resolutions \cite{BeJo}. One natural question to ask is whether the Tate-Hochschild cohomology shares the same properties as those of the usual Hochschild cohomology. In the following sections, we will examine the Tate-Hochschild cohomology of $RG$ and its ring structure.

\section{Tate-Hochschild Cohomology of a Group Algebra}
\label{sec:Tate-Hochschild RG}

The opposite algebra of $RG$, denoted by $RG^{op}$, is $RG$ with the same addition but with multiplication performed in the reverse order. Let $RG^e = RG \otimes RG^{op}$ be the enveloping algebra of $RG$. Since $G$ is finite, $RG^e$ is two-sided Noetherian and Gorenstein of Gorenstein dimension 0. Let $M$ be an $RG$-bimodule. $M$ may be regarded as a left $RG^e$-module by setting $(a \otimes b)\cdot m = amb$, where $a \in RG, b \in RG^{op}$, and $m \in M$. In particular, we shall regard $RG$ as a left $RG^e$-module. For any integer $n \in \mathbb{Z}$, the $n$-th Tate-Hochschild cohomology of $RG$ is defined as:
 \[\widehat{\HH}^n(RG,M) := \widehat{\Ext}^n_{RG^e}(RG,M), \]
where the $\widehat{\Ext}$ functor is taken using an $RG^e$-complete resolution of $RG$. As $RG^e$ is a two-sided Noetherian and Gorenstein ring, Theorems 3.1 and 3.2 in \cite{AvMa} guarantee that every finitely generated $RG^e$-module admits a complete resolution. Hence, we obtain an $RG^e$-complete resolution $\mathbb{X}$ for $RG$. The $n$-th Tate-Hochschild cohomology group of $RG$ is the $n$-th homology group of the complex $\Hom_{RG^e}(\mathbb{X},M)$. 

\begin{Remark}{\em 
The Tate-Hochschild cohomology groups of $RG$ agree with the usual Hochschild cohomology groups in all positive degrees:
 \[\widehat{\HH}^n(RG,M) \cong \HH^n(RG,M), \text{ for all } n>0. \]}
\end{Remark}

Let $M$ and $N$ be $RG$-bimodules. Then $M \otimes_{RG} N$ is also an $RG$-bimodule, which can be considered as a left $RG^e$-module via $(a \otimes b) \cdot (m \otimes_{RG} n) = am \otimes_{RG} nb$, where $a \in RG, b \in RG^{op}, m \in M$, and $n \in N$. There is a cup product on Tate-Hochschild cohomology:
 \[\widehat{\HH}^i(RG,M) \otimes \widehat{\HH}^j(RG,N) \rightarrow \widehat{\HH}^{i+j}(RG,M \otimes_{RG} N), \]
which we describe as follows. Let $\mathbb{X}$ be any $RG^e$-complete resolution of $RG$. The complete tensor product $\mathbb{X} \widehat{\otimes}_{RG} \mathbb{X}$ is obtained by defining:
 \[(\mathbb{X} \widehat{\otimes}_{RG} \mathbb{X})_n = \displaystyle \prod_{i+j=n} X_i \otimes_{RG} X_j, \text{ for all } n \in \mathbb{Z}. \]
By the same argument as in (\cite{Brown}, Section~VI.5) or (\cite{Nguyen}, Section~6.1), $\mathbb{X} \widehat{\otimes}_{RG} \mathbb{X}$ is an acyclic chain complex of $RG^e$-modules. Lemma~6.5 in \cite{Nguyen} shows the existence of a complete diagonal approximation chain map $\Gamma: \mathbb{X} \rightarrow \mathbb{X} \widehat{\otimes}_{RG} \mathbb{X}$ that preserves the augmentation.

Given graded modules $B, B', C, C'$ and module homomorphisms $u: C \rightarrow B$ of degree $i$ and $v: C' \rightarrow B'$ of degree $j$, there is a map $u \widehat{\otimes} v: C \widehat{\otimes} C' \rightarrow B \widehat{\otimes} B'$ of degree $i+j$ defined by:
 \[ (u \widehat{\otimes} v)_n = \displaystyle \prod_{r+s=n} (-1)^{rj} u_r \otimes v_s: \displaystyle \prod_{r+s=n} C_r \otimes C'_s \rightarrow \displaystyle \prod_{r+s=n} B_{r+i} \otimes B'_{s+j}. \]  
 
Let $f \in \Hom_{RG^e}(X_i,M)$ represent an element of $\widehat{\HH}^i(RG,M)$ and let $g \in \Hom_{RG^e}(X_j,N)$ represent an element of $\widehat{\HH}^j(RG,N)$. Then:
 \[ f \smile g = (f \widehat{\otimes} g) \circ \Gamma \in \Hom_{RG^e}(\mathbb{X},M \otimes_{RG} N) \] 
represents an element of $\widehat{\HH}^{i+j}(RG,M \otimes_{RG} N)$. One can check that this product is independent of $\mathbb{X}$ and $\Gamma$ and satisfies the usual cup product properties. When $M=N=RG$, this cup product gives $\widehat{\HH}^*(RG,RG)$ the structure of an associative graded ring, as $RG \cong RG \otimes_{RG} RG$. 

\section{Reduction to Tate Cohomology and Relations with Subgroups}
\label{sec:relationship RG}

We show here that the Tate-Hochschild cohomology of $RG$ can be reduced to its Tate cohomology. We begin with a lemma which is based on the original Eckmann-Shapiro Lemma but is generalized to a complete resolution:

\begin{lemma}[Eckmann-Shapiro, (\cite{Nguyen}, Lemma~7.1) or (\cite{Brown}, VI.5.2)]
\label{Eckmann-Shapiro RG}
Let $H$ be a subgroup of a finite group $G$, $M$ be a left $RH$-module and $N$ be a left $RG$-module. Consider $N\downarrow^G_H=N$ to be a left $RH$-module via restriction of the action, and let $M\uparrow^G_H:= RG \otimes_{RH} M$ denote the induced $RG$-module where $G$ acts on the leftmost factor by multiplication. Then for all $n \in \mathbb{Z}$, there is an isomorphism of abelian groups:
 \[\widehat{\Ext}^n_{RH}(M,N\downarrow^G_H) \cong \widehat{\Ext}^n_{RG}(M\uparrow^G_H,N). \]
\end{lemma}

For any subgroup $H \subseteq G$, every $RH$-complete resolution gives rise to an $RG$-complete resolution by inducing the modules (induction takes projectives to projectives, and exact sequences to exact sequences). Since coinduction is the same as induction in the finite group case (\cite{Brown}, Prop.~III.5.9), the proof for this lemma goes through in the present context and follows from the Nakayama relations (\cite{Ben1}, Prop.~2.8.3). 

We have $RG \cong RG^{op}$ as algebras via $g \mapsto g^{-1}$, and $RG \otimes RG \cong R(G \times G)$ as algebras. As a result, $RG^e = RG \otimes RG^{op} \cong R(G \times G)$. Here, $RG$ is a left $R(G \times G)$-module by the two-sided action: $(g_1,g_2) \cdot g = g_1gg_2^{-1}$. So $RG$ is just the permutation module $R\uparrow^{G \times G}_{\delta(G)} = R(G \times G) \otimes_{R\delta(G)} R$ on the cosets of the diagonal $\delta(G)=\{(g,g)\:|\:g \in G\}$. Moreover, any $RG$-bimodule $M$ may be considered as a left $RG$-module with action given by $g \cdot m = gmg^{-1}$, for all $g \in G$ and $m \in M$. It follows from this discussion and Lemma \ref{Eckmann-Shapiro RG} that:
\begin{align*}
\widehat{\HH}^*(RG,M) = \widehat{\Ext}^*_{RG^e}(RG,M) &\cong \widehat{\Ext}^*_{R(G \times G)}(RG,M) \\
&\cong \widehat{\Ext}^*_{R(G \times G)}(R\uparrow^{G \times G}_{\delta(G)},M) \\
&\cong \widehat{\Ext}^*_{RG}(R,M\downarrow^{G \times G}_{\delta(G)}) \\
&\cong \widehat{\Ho}^*(G,M).
\end{align*}
The Tate-Hochschild cohomology of $G$ with coefficients in the bimodule $M$ is just the Tate cohomology of $G$ with coefficients in $M$ under the diagonal action. When $M$ is an algebra itself with compatible $G$-action, this is in fact an algebra isomorphism as it respects the cup products defined for Tate and Tate-Hochschild cohomology (\cite{Nguyen}, Theorem~7.2). In particular, when $M=RG$ with $G$-action via conjugation,
\begin{equation}
\label{THdecomp eq:1}
\widehat{\HH}^*(RG,RG) \cong \widehat{\Ho}^*(G,RG). 
\end{equation}
Therefore, all properties of $\widehat{\Ho}^*(G,RG)$ that we observed in Section~\ref{sec:Tate RG} transfer to those for the Tate-Hochschild cohomology. 

Let $H$ be a subgroup of $G$. By restricting the action, any $RG$-module $N$ may be regarded as an $RH$-module and any $RG$-complete resolution $\mathbb{P}$ of $R$ may also be considered as an $RH$-complete resolution. Sections~XII.8 and XII.9 in \cite{CaEi} show that there are maps in the Tate cohomology with properties analogous to those in the usual group cohomology:
\begin{itemize}
 \item The restriction map:
 \[\res^G_H: \widehat{\Ho}^*(G,N) \rightarrow \widehat{\Ho}^*(H,N),\]
which is induced from the inclusion $\Hom_{RG}(\mathbb{P},N) \subset \Hom_{RH}(\mathbb{P},N)$.
 \item The corestriction map (or transfer):
 \[\cor^G_H: \widehat{\Ho}^*(H,N) \rightarrow \widehat{\Ho}^*(G,N),\]
which is given on the cochain level by defining: 
 \[(\cor^G_H f)(p)=\displaystyle \sum_{g\in \mathcal{G}}gf(g^{-1}p),\]
where $\mathcal{G}$ denotes a set of left coset representatives of $H$ in $G$, $f \in \Hom_{RH}(P_i,N)$, and $p \in P_i$. One can check that this definition is independent of the choice of the representatives $g \in \mathcal{G}$.
 \item Moreover, for any $g \in G$, there is an isomorphism:
 \[g^*: \widehat{\Ho}^*(H,N) \rightarrow \widehat{\Ho}^*(gHg^{-1},N) = \widehat{\Ho}^*(^gH,N)\]
defined on the cochain level as $(g^*f)(p)=g(f(g^{-1}p)$.
\end{itemize}
We shall list some properties of these maps without proving them. The proof goes through using similar arguments as in \cite{CaEi}. The readers can refer to \cite{CaEi} for more details.

\begin{proposition}[\cite{CaEi}, Section~XII.8 (4)-(14) and Section~XII.9 (4)]
\label{THdecomp prop}
Let $K \subseteq H \subseteq G$ be subgroups, and $N_1, N_2$ be $RG$-modules which may be regarded as $RH$-modules. Let $\alpha_i \in \widehat{\Ho}^*(G,N_i)$, $\beta_i \in \widehat{\Ho}^*(H,N_i)$, and $g_i \in G$, for $i=1,2$. Then the maps defined above satisfy:
\begin{enumerate}
 \item $g_1^*g_2^*=(g_1g_2)^*$
 \item $g^* = \textbf{1} \text{, if } g \in H$
 \item $\cor^G_H \circ \res^G_H = (G:H)\textbf{1}$
 \item $\res^H_K \circ \res^G_H = \res^G_K$
 \item $\cor^G_H \circ \cor^H_K = \cor^G_K$
 \item $g^* \circ \res^H_K = \res^{^gH}_{^gK} \circ g^*$
 \item $g^* \circ \cor^H_K = \cor^{^gH}_{^gK} \circ g^*$
 \item $\res^G_H(\alpha_1 \smile \alpha_2)= (\res^G_H \alpha_1) \smile (\res^G_H \alpha_2)$
 \item $\cor^G_H (\beta_1 \smile \res^G_H \alpha_2) = (\cor^G_H \beta_1) \smile \alpha_2$
 \item $\cor^G_H (\res^G_H \alpha_1 \smile \beta_2) = \alpha_1 \smile (\cor^G_H \beta_2)$
 \item $g^*(\beta_1 \smile \beta_2) = (g^*\beta_1) \smile (g^*\beta_2)$
 \item Let $H, K \subseteq G$ be subgroups and $N$ be an $RG$-module which may be regarded as an $RH$ (or $RK$)-module. The map $\res^G_K \circ \cor^G_H: \widehat{\Ho}^*(H,N) \rightarrow \widehat{\Ho}^*(K,N)$ is given by: 
 \[\res^G_K (\cor^G_H (\beta))= \displaystyle \sum_{x \in D} \cor^K_{K \cap \,^xH}(\res^{^xH}_{K \cap \,^xH}(x^*\beta)),\]
 where $\beta \in \widehat{\Ho}^*(H,N)$ and $D$ is a set of double coset representatives such that $G=\displaystyle \bigcup_{x \in D} KxH$ is a disjoint union.
\end{enumerate}
\end{proposition}

\section{Generalized Additive Decomposition}
\label{sec:decomp}

For the rest of this paper, we let $H$ be another finite group which acts as automorphisms on $G$. Via this action, $RG$ becomes an $RH$-module. The multiplication map $RG \otimes RG \rightarrow RG$ is an $RH$-module homomorphism. Hence, it induces the ring structure on cohomology $\widehat{\Ho}^*(H,RG):= \widehat{\Ext}^*_{RH}(R,RG)$ by compositing with the cup product. We will study the additive decomposition of this ring $\widehat{\Ho}^*(H,RG)$. The Tate-Hochschild cohomology ring $\widehat{\Ho}^*(RG,RG) \cong \widehat{\Ho}^*(G,RG)$ is a special case of $\widehat{\Ho}^*(H,RG)$ by letting $H=G$ act on itself by conjugation. 

\begin{proposition}
\label{abelian RG}
If $H$ acts trivially on $G$, then $\widehat{\Ho}^*(H,RG) \cong RG \otimes_R \widehat{\Ho}^*(H,R)$ as graded $R$-algebras. In particular, if $G$ is abelian, then   \[\widehat{\HH}^*(RG,RG) \cong RG \otimes_R \widehat{\Ho}^*(G,R). \]
\end{proposition}

\vspace{-1em}
\begin{proof}
If $G$ is abelian, $H=G$ acting on itself by conjugation yields the trivial action. Hence, the second statement follows from the first statement and the isomorphism (\ref{THdecomp eq:1}):
 \[\widehat{\HH}^*(RG,RG) \cong \widehat{\Ho}^*(G,RG) \cong RG \otimes_R \widehat{\Ho}^*(G,R). \]
 
To prove the first statement, let $\varepsilon: \mathbb{P} \rightarrow R$ be an $RH$-complete resolution of $R$. Since $H$ acts trivially on $G$, $RG$ is a trivial $RH$-module and is free as an $R$-module.  We claim that $\gamma: RG \otimes \Hom_{RH}(\mathbb{P},R) \rightarrow \Hom_{RH}(\mathbb{P},RG)$ is an isomorphism. It can be seen by sending $g \otimes f \mapsto \gamma(g \otimes f)=F$, where $F(p)=f(p)g$, for $f \in \Hom_{RH}(P_i,R), p \in P_i$, and $g \in G$. It is easy to check that $F \in \Hom_{RH}(\mathbb{P},RG)$ and it is a cocycle when $f$ is. Hence, passing to the homology, $\gamma$ induces an isomorphism of graded $R$-modules:
 \[\gamma_*: RG \otimes_R \widehat{\Ho}^*(H,R) \rightarrow \widehat{\Ho}^*(H,RG). \]
The definition of cup product corresponds to this map, making $\gamma_*$ a ring homomorphism. 
\end{proof}

This proposition helps us to find and study the structure of the Tate-Hochschild cohomology ring of a finite abelian group algebra, given its Tate cohomology ring. For example, knowing the Tate cohomology of a cyclic group $G$, (see \cite{CaEi}, Section~XII.7), one can easily compute its Tate-Hochschild cohomology by applying Proposition~\ref{abelian RG}. 

Now we return to the general case where $H$ acts on $G$ non-trivially and $G$ is not necessarily abelian. Let $g_1, \ldots, g_t \in G$ be representatives of the orbits of the action of $H$ on $G$. Let $H_i:=\text{Stab}_H(g_i)=\{h \in H \:|\: ^hg_i=g_i\}$ be the stabilizer of $g_i$. For any $g \in G$, there are two $R(\text{Stab}_H(g))$-module homomorphisms:
\[\theta_g: R \rightarrow RG \text{ via } r \mapsto rg, \]
\[\pi_g: RG \rightarrow R \text{ via } \sum_{a\in G} r_aa \mapsto r_g. \]
If $V$ is any subgroup of $\text{Stab}_H(g)$, then these maps induce maps on cohomology:
\[\theta_g^*: \widehat{\Ho}^*(V,R) \rightarrow \widehat{\Ho}^*(V,RG), \]
\[\pi_g^*: \widehat{\Ho}^*(V,RG) \rightarrow \widehat{\Ho}^*(V,R), \]
since $\widehat{\Ext}^*$ is covariant in the second argument. The following properties of $\theta_g^*$ and $\pi_g^*$ will help us in proving the main result. 

\begin{lemma}
\label{THdecomp lemma}
Let $h \in H$ and $a, b \in G$.
\begin{enumerate}
 \renewcommand{\labelenumi}{(\alph{enumi})} 
 \item If $V$ is a subgroup of $\text{Stab}_H(a)$, then $h^* \circ \theta_a^* = \theta_{^ha}^* \circ h^*$ as maps from $\widehat{\Ho}^*(V)$ to $\widehat{\Ho}^*(^hV,RG)$.
 \item Suppose $V \subseteq \text{Stab}_H(a) \cap \text{Stab}_H(b)$ and $\alpha, \beta \in \widehat{\Ho}^*(V)$. Then:
 \[\theta_a^*(\alpha) \smile \theta_b^*(\beta) = \theta_{ab}^*(\alpha \smile \beta). \]
 \item Suppose $V' \subseteq V \subseteq \text{Stab}_H(a)$. Then $\theta_a^*$ and $\pi_a^*$ commute with $\res^V_{V'}$ and $\cor^V_{V'}$.
 \item If $V \subseteq \text{Stab}_H(a) \cap \text{Stab}_H(b)$, then $\pi_a^* \circ \theta_b^* = \delta_{a,b}\textbf{1}$, where $\textbf{1}$ is the identity map on $\widehat{\Ho}^*(V)$ and $\delta_{a,b}$ is the Kronecker delta.
\end{enumerate}
\end{lemma}

\vspace{-1em}
\begin{proof}
Lemma~5.2 in \cite{SW} showed these properties in positive degrees. We generalize the proof to negative degrees and present it on the cochain level. The desired results are induced on cohomology.
\begin{enumerate}
 \renewcommand{\labelenumi}{(\alph{enumi})} 
 \item Let $\mathbb{P}$ be an $RV$-complete resolution of $R$, $f \in \Hom_V(P_i,R)$ be a cocycle representing an element of $\widehat{\Ho}^i(V)$, and $p \in P_i$. Then
 \[h^*(\theta_a f)(p)= f(h^{-1}p)(^ha)= \theta_{^ha}(h^*(f))(p). \]
 
 \item Let $m: RG \otimes RG \rightarrow RG$ be the multiplication map and $\Gamma: \mathbb{P} \rightarrow \mathbb{P} \widehat{\otimes} \mathbb{P}$ be a complete diagonal approximation map. Let $f, q \in \Hom_{V}(\mathbb{P},R)$ represent $\alpha, \beta \in \widehat{\Ho}^*(V)$, respectively. Then on the cochain level:
 \[m \circ ((\theta_a \circ f) \widehat{\otimes} (\theta_b \circ q)) \circ \Gamma = m \circ (\theta_a \otimes \theta_b) \circ (f \widehat{\otimes} q) \circ \Gamma = \theta_{ab} \circ (f \widehat{\otimes} q) \circ \Gamma, \]
where the left side represents $\theta_a^*(\alpha) \smile \theta_b^*(\beta)$ and the right side represents $\theta_{ab}^*(\alpha \smile \beta)$.

 \item Let $\mathbb{P}$ be an $RV$-complete resolution of $R$ which can also be regarded as an $RV'$-complete resolution by restricting the action. Let $f \in \Hom_{V'}(P_i,RG)$ represent an element of $\widehat{\Ho}^i(V',RG)$, $q \in \Hom_{V'}(P_i,R)$ represent an element of $\widehat{\Ho}^i(V')$, and $p \in P_i$.
 \begin{align*}\displaystyle (\pi_a^* \cor^V_{V'})(f)(p)= \pi_a\left(\sum_{v \in V/V'} (^vf(v^{-1}p))\right) &= \sum_v \pi_{^{v^{-1}}a} (f(v^{-1}p)) \\
 &= \sum_v (\pi_a \circ f)(v^{-1}p)= (\cor^V_{V'} \pi_a^*)(f)(p),
 \end{align*}
since $v \in V \subseteq \text{Stab}_H(a)$, we have $^{v^{-1}}a=a$, and $V$ acts trivially on $R$. Similarly,
 \begin{align*}\displaystyle (\theta_a^* \cor^V_{V'})(q)(p)= \theta_a\left(\sum_{v \in V/V'} q(v^{-1}p)\right) &= \sum_v \theta_{^{v^{-1}}a} (q(v^{-1}p))\\
 &= \sum_v (\theta_a \circ q)(v^{-1}p)= (\cor^V_{V'} \theta_a^*)(q)(p). 
 \end{align*}
The other cases follow similarly by commutativity between $\pi_a$, $\theta_a$ and the inclusion map $\iota: \Hom_V(\mathbb{P},N) \hookrightarrow \Hom_{V'}(\mathbb{P},N)$, where $N=RG$ or $R$. 

 \item Let $r \in R$.
 $$\pi_a(\theta_b(r))= \pi_a(rb)= \begin{cases}
                                  r, &\quad \text{if } a=b \\
                                  0, &\quad \text{else}.
                                  \end{cases}$$
\end{enumerate}
\end{proof}
\vspace{-1em}

For all $i=1,2,\ldots,t$, let $\psi_i: \widehat{\Ho}^*(H_i,R) \rightarrow \widehat{\Ho}^*(H,RG)$ be defined as the composition $\psi_i= \cor^H_{H_i} \circ \theta_{g_i}^*$. We describe the additive decomposition that is generalized from the usual cohomology (\cite{Ben2}, Theorem~2.11.2):

\begin{lemma}
\label{THdecomp}
The map $\displaystyle \widehat{\Ho}^*(H,RG) \rightarrow \bigoplus_i \widehat{\Ho}^*(H_i,R)$, sending $\zeta \mapsto (\pi_{g_i}^* \circ \res^H_{H_i}(\zeta))_i$, is an isomorphism of graded $R$-modules, for $\zeta \in \widehat{\Ho}^*(H,RG)$. Its inverse sends $\alpha \in \widehat{\Ho}^*(H_i,R)$ to $\psi_i(\alpha) \in \widehat{\Ho}^*(H,RG)$.
\end{lemma}

\vspace{-1em}
\begin{proof}
For $i=1,2,\ldots,t$, let $M_i$ be the free $R$-module generated by elements of the orbit containing $g_i$. Then $RG= \bigoplus_i M_i$. There is an isomorphism $M_i \rightarrow R\uparrow^H_{H_i}=RH \otimes_{RH_i} R$ given by $r(^hg_i) \mapsto h \otimes r$. It induces an isomorphism in cohomology $\widehat{\Ho}^*(H,M_i) \cong \widehat{\Ho}^*(H,R\uparrow^H_{H_i})$.

Since $\widehat{\Ext}$ is additive, $\widehat{\Ho}^*(H,RG) \cong \bigoplus_i \widehat{\Ho}^*(H,M_i)$. Apply the Eckmann-Shapiro Lemma \ref{Eckmann-Shapiro RG}, we have $\widehat{\Ho}^*(H,RG) \cong \bigoplus_i \widehat{\Ho}^*(H_i,R)$. One can also check directly that the maps given in the statement of the lemma are inverses of each other by taking their compositions and applying Proposition~\ref{THdecomp prop} and Lemma~\ref{THdecomp lemma} to show that their compositions are the identity maps. 
\end{proof}

\begin{Remark}{\em
If $H=G$ acts on itself by conjugation, then $M_i$ is the free $R$-module generated by the conjugacy class of $g_i$. $M_i$ is isomorphic to $R\uparrow^G_{C_G(g_i)}$, where $C_G(g_i)$ is the centralizer of $g_i$. Therefore, the isomorphism in Lemma~\ref{THdecomp} gives an additive decomposition of the Tate-Hochschild cohomology of $G$ as a direct sum of the Tate cohomology of the centralizers of conjugacy class representatives of $G$ with coefficients in $R$:
\[ \widehat{\HH}^*(RG,RG) \cong \bigoplus_i \widehat{\Ho}^*(C_G(g_i),R). \] }
\end{Remark}
\vspace{-0.4em}
In 1999, Siegel and Witherspoon showed that there is a product formula for the usual Hochschild cohomology of $G$ in terms of a similar additive decomposition (\cite{SW}, Theorem~5.1). We will describe products in $\widehat{\Ho}^*(H,RG)$ with respect to the isomorphism in Lemma~\ref{THdecomp}. The argument will be analogous to that in \cite{SW} but is generalized to the Tate cohomology.

Fix $i, j \in \{1,2,\ldots,t\}$. Let $D$ be a set of double coset representatives for $H_i\backslash H/H_j$. For each $x \in D$, there is a unique $k=k(x)$ such that
\begin{equation}
\label{THdecomp eq:2}
g_k=\, ^yg_i\,^{yx}g_j
\end{equation} 
for some $y \in H$. One can expand the action on the right hand side and get $g_k=\, ^y(g_i\,^xg_j)$ showing that $g_k$ is just a representative of an orbit of the action of $H$ on the double coset $H_ixH_j$. Moreover, $k$ is independent of the choice of double coset representative $x$. The set of all $y$ satisfying (\ref{THdecomp eq:2}) is also a double coset. To see this, let us fix $y=y(x)$ for which (\ref{THdecomp eq:2}) holds. Let $y' \in H$ be another element such that $g_k=\, ^{y'}g_i\,^{y'x}g_j$. Then $^{y'}g_i\, ^{y'x}g_j= g_k=\, ^yg_i\, ^{yx}g_j$ implies $^{y'}(g_i\,^xg_j)=\, ^y(g_i\,^xg_j)$. Let $h=y'y^{-1}$. We have $^hg_k=\, ^h(^y(g_i\,^xg_j))= \,^{y'}(g_i\,^xg_j)= g_k$ showing $h \in H_k = \text{Stab}_H(g_k)$. On the other hand, if $h \in H_k$, then let $y'=hy \in H$, we have $^{hy}g_i\, ^{hyx}g_j= \,^h(^yg_i\, ^{yx}g_j) = \,^hg_k = g_k$. Putting together, we have shown: 
\[\{ y' \in H|g_k= \,^{y'}g_i\, ^{y'x}g_j\} = H_ky = H_ky(^xH_j \cap H_i) \in H_k\backslash H/(^xH_j \cap H_i), \]
where the last equality follows from (\ref{THdecomp eq:2}) and $^{yx}H_j \cap \,^yH_i \subseteq H_k$. We can now prove our main result which provides a formula for products in $\widehat{\Ho}^*(H,RG)$ with respect to Lemma~\ref{THdecomp}.

\begin{theorem}
\label{THdecomp product}
Let $\alpha \in \widehat{\Ho}^*(H_i)$ and $\beta \in \widehat{\Ho}^*(H_j)$. Then
\[\psi_i(\alpha) \smile \psi_j(\beta) = \sum_{x \in D} \psi_k(\cor^{H_k}_V(\res^{^yH_i}_V y^*\alpha \smile \res^{^{yx}H_j}_V (yx)^*\beta)) \]
where $D$ is a set of double coset representatives for $H_i \backslash H/H_j$, $k=k(x)$ and $y=y(x)$ are chosen to satisfy (\ref{THdecomp eq:2}), and $V=V(x)= \,^{yx}H_j \cap \,^yH_i \subseteq H_k$.
\end{theorem}

\vspace{-1em}
\begin{proof}
By Lemma~\ref{THdecomp},
\begin{align*}
\psi_i(\alpha) \smile \psi_j(\beta) &= \cor^H_{H_i}(\theta_{g_i}^* \alpha) \smile \cor^H_{H_j}(\theta_{g_j}^* \beta), \text{ by definition of } \psi_i, \psi_j \\
&= \cor^H_{H_i}(\theta_{g_i}^* \alpha \smile \res^H_{H_i}\cor^H_{H_j}\theta_{g_j}^* \beta), \text{ by Prop.~\ref{THdecomp prop} (9)} \\
&= \sum_{x \in D} \cor^H_{H_i}(\theta_{g_i}^* \alpha \smile \cor^{H_i}_{^xH_j \cap H_i}\res^{^xH_j}_{^xH_j \cap H_i}x^*\theta_{g_j}^* \beta), \text{ by Prop.~\ref{THdecomp prop} (12)} \\
&= \sum_{x \in D} \cor^H_{H_i}(\cor^{H_i}_{^xH_j \cap H_i}(\res^{H_i}_{^xH_j \cap H_i}\theta_{g_i}^* \alpha \smile \res^{^xH_j}_{^xH_j \cap H_i}x^*\theta_{g_j}^* \beta)), \text{ by Prop.~\ref{THdecomp prop} (10)} \\
&= \sum_{x \in D} \cor^H_{^xH_j \cap H_i}(\res^{H_i}_{^xH_j \cap H_i}\theta_{g_i}^* \alpha \smile \res^{^xH_j}_{^xH_j \cap H_i}x^*\theta_{g_j}^* \beta), \text{ by Prop.~\ref{THdecomp prop} (5)} \\
&= \sum_{x \in D} \cor^H_{^xH_j \cap H_i} \theta_{g_i\,^xg_j}^*(\res^{H_i}_{^xH_j \cap H_i} \alpha \smile \res^{^xH_j}_{^xH_j \cap H_i}x^*\beta), \text{ by Lemma~\ref{THdecomp lemma} (a)-(c)} \\
&= \sum_k \sum_x \psi_k \pi_{g_k}^* \res^H_{H_k} (\cor^H_{^xH_j \cap H_i} \theta_{g_i\,^xg_j}^*(\res^{H_i}_{^xH_j \cap H_i} \alpha \smile \res^{^xH_j}_{^xH_j \cap H_i} x^*\beta)), \\
& \text{ by the isomorphism in Lemma~\ref{THdecomp}} \\
&= \sum_k \sum_{x,y} \psi_k \pi_{g_k}^* \cor^{H_k}_{V'} \res^{^{yx}H_j \cap \,^yH_i}_{V'} y^* \theta_{g_i\,^xg_j}^*(\res^{H_i}_{^xH_j \cap H_i} \alpha \smile \res^{^xH_j}_{^xH_j \cap H_i}x^*\beta), \\
& \text{ by Prop.~\ref{THdecomp prop} (12), where $y$ runs over a set of representatives for $H_k\backslash H/^xH_j \cap H_i$} \\
& \text{ and $V'=H_k \cap \,^{yx}H_j \cap \,^yH_i$,}\\
&= \sum_k \sum_{x,y} \psi_k \cor^{H_k}_{V'} \pi_{g_k}^* \theta_{^yg_i\,^{yx}g_j}^* \res^{^{yx}H_j \cap \,^yH_i}_{V'} y^*(\res^{H_i}_{^xH_j \cap H_i} \alpha \smile \res^{^xH_j}_{^xH_j \cap H_i}x^*\beta), \\
& \text{ by Lemma~\ref{THdecomp lemma} (a),(c)} \\
&= \sum_{x \in D} \psi_k \cor^{H_k}_{V'} (\res^{^yH_i}_{V'} y^* \alpha \smile \res^{^{yx}H_j}_{V'} (yx)^* \beta), \\
& \text{ by Prop.~\ref{THdecomp prop} (1), (4), (6) and Lemma~\ref{THdecomp lemma} (d)}.
\end{align*} 
By Lemma~\ref{THdecomp lemma} (d), the only terms that can be non-zero in the next to last step are those for which $g_k= \,^yg_i\, ^{yx}g_j$. We have seen in the discussion prior to this theorem that each $x$ determines a unique $k$ and double coset $H_ky(^xH_j \cap H_i)$ for which this holds. Therefore, we may take $y=y(x)$ and  $^{yx}H_j \cap \,^yH_i \subseteq H_k$. Hence, $V'=V= \,^{yx}H_j \cap \,^yH_i$.
\end{proof}

\begin{Remark}{\em
Since the cup product is well-defined and unique (\cite{CaEi}, Theorem~XII.5.1), the sum in the statement of the theorem is independent of the choice of $x$ and $y$. One can see this directly by replacing $y$ with $hy$, for some $h \in H_k$. By Prop.~\ref{THdecomp prop} (6), (7), and (11), $h^*$ respects the cup product and commutes with the restriction and corestriction maps. Moreover, since $H_k$ acts trivially on its own cohomology, any term of the sum in the theorem is unchanged by replacing $y$ with $hy$. If $x$ is multiplied on the right by an element of $H_j$, the terms are unchanged for similar reasons. If $x$ is replaced by $hx$, for some $h \in H_i$, then we must replace $y$ with $yh^{-1}$ so that (\ref{THdecomp eq:2}) holds:
\[^{(yh^{-1})}g_i\, ^{(yh^{-1})(hx)}g_j = \,^yg_i\, ^{yx}g_j = g_k, \]
and the terms remain unchanged. }
\end{Remark}
We observe that when $i=1$, $\psi_1: \widehat{\Ho}^*(H,R) \rightarrow \widehat{\Ho}^*(H,RG)$ is an algebra monomorphism that is induced by the algebra homomorphism $R \rightarrow RG$ mapping $r \mapsto r1$. Alternatively, by letting $i=j=1$ in Theorem \ref{THdecomp product}, we see that $\psi_1$ respects the cup product:
\[\psi_1(\alpha) \smile \psi_1(\beta)= \psi_1(\alpha \smile \beta), \]
where $\alpha, \beta \in \widehat{\Ho}^*(H)$. Hence, via $\psi_1$, we may view $\widehat{\Ho}^*(H,RG)$ as a (left) $\widehat{\Ho}^*(H)$-module with action via multiplying (on the left) by $\psi_1(\widehat{\Ho}^*(H))$. Each $\widehat{\Ho}^*(H_i)$ may also be regarded as an $\widehat{\Ho}^*(H)$-module via restriction. As a consequence, we obtain:

\begin{corollary}
The isomorphism in Lemma \ref{THdecomp} is an isomorphism of graded $\widehat{\Ho}^*(H)$-modules:
\[\widehat{\Ho}^*(H,RG) \stackrel{\cong}{\longrightarrow} \bigoplus_i \widehat{\Ho}^*(H_i) \] 
\end{corollary}

\vspace{-1em}
\begin{proof}
For $i=1$, let  $\alpha \in \widehat{\Ho}^*(H)$ and $\beta \in \widehat{\Ho}^*(H_j)$. Theorem \ref{THdecomp product} reduces to:
\[\psi_1(\alpha) \smile \psi_j(\beta)= \psi_j(\res^H_{H_j} (\alpha) \smile \beta). \]
where the left hand side is considered as action of $\widehat{\Ho}^*(H)$ on $\widehat{\Ho}^*(H,RG)$ that corresponds to the action of $\widehat{\Ho}^*(H)$ on each $\widehat{\Ho}^*(H_j)$ on the right hand side, via the isomorphism in Lemma \ref{THdecomp}.
\end{proof}

As noted in the remark following Lemma \ref{THdecomp}, when $H=G$ acts on itself by conjugation, Theorem \ref{THdecomp product} gives a formula for the multiplicative structure of $\widehat{\HH}^*(RG,RG) \cong \bigoplus_i \widehat{\Ho}^*(C_G(g_i),R)$ in terms of this decomposition. It reduces the computation of products in $\widehat{\HH}^*(RG,RG)$ to products within the Tate cohomology rings of certain subgroups of $G$. In the next section, we will show a basic non-abelian example that demonstrates the usefulness of Theorem \ref{THdecomp product}.

\section{The symmetric group on three elements}
\label{sec:S3}

Let $\bk$ be a field of characteristic 3. Let $G=S_3=\left\langle a, b \: |\: a^3=1=b^2, ab=ba^2 \right\rangle$ act on itself by conjugation. Without loss of generality, we choose conjugacy class representatives $g_1=1$, $g_2=a$, and $g_3=b$ whose centralizers are $H_1=G$, $H_2=\{1, a, a^2\}=:N$, and $H_3=\{1,b\}$, respectively. We will find the Tate-Hochschild cohomology ring of $\kG$ using elements of $\widehat{\Ho}^*(H_i,\bk)$ and the product formula given in Theorem \ref{THdecomp product}. 

Let us examine each ring $\widehat{\Ho}^*(H_i):= \widehat{\Ho}^*(H_i,\bk)$. Since the characteristic of $\bk$ is $3$ and $N$ is cyclic of order $3$, the cohomology ring $\widehat{\Ho}^*(N)$ is periodic by (\cite{CaEi}, Theorem~XII.11.6) and (\cite{Ben}, (4.1.3)). By direct computation from Section~XII.7 in \cite{CaEi}, $\widehat{\Ho}^*(N)$ is of the form $\Lambda(w_1) \otimes_{\bk} \bk[w_2,w_2^{-1}]$, where $\Lambda(w_1)$ is the exterior $\bk$-algebra on the element $w_1$ of degree $1$ and $\bk[w_2,w_2^{-1}]$ is generated by the elements $w_2$ of degree $2$ and $w_2^{-1}$ of degree $-2$, subject to the graded-commutative relations and $w_2w_2^{-1}=1$. By a similar computation, because the characteristic of $\bk$ does not divide the order of $H_3=\{1,b\}$, we find that $\widehat{\Ho}^n(H_3)=0$ in all degrees. As a result, $\widehat{\Ho}^*(H_3)=0$. 

We now compute the Tate cohomology ring of $G=S_3$ with coefficients in $\bk$. $N$ is a normal subgroup of $G$. The quotient group $G/N$ is (isomorphic to) $\mathbb{Z}_2$. It is easy to check that $G$ is isomorphic to a semidirect product $N \rtimes \mathbb{Z}_2$ and every abelian subgroup of $G$ is cyclic. It follows from (\cite{CaEi}, Theorem~XII.11.6) and (\cite{Ben}, (4.1.3)) that the Tate cohomology ring $\widehat{\Ho}^*(G)$ is periodic and Noetherian. One can directly compute $\widehat{\Ho}^*(G)$ by using an $N$-complete resolution of $\bk$, imposing on it an action of $\mathbb{Z}_2$ to make it become a $G$-complete resolution of $\bk$, computing the Tate cohomology groups from that resolution, and studying their products. Alternatively, following the discussion in (\cite{CaEi}, Section~XII.10), we see that for any $G$-module $M$, $\widehat{\Ho}^*(G,M)$ is a direct sum of $\widehat{\Ho}^*(G,M,p)$, where $\widehat{\Ho}^*(G,M,p)$ is the $p$-primary component of $\widehat{\Ho}^*(G,M)$ and $p$ runs through all the prime divisors of $|G|=6$. Here, $M=\bk$ is a field of characteristic $3$, so only the $3$-primary component is non-zero. By (\cite{CaEi}, Theorem~XII.10.1), $G/N$ operates on $\widehat{\Ho}^*(N)$ and so $\widehat{\Ho}^*(G)=\widehat{\Ho}^*(G,\bk,3) \cong \left[\widehat{\Ho}^*(N)\right]^{G/N} \cong \Lambda(w_1w_2) \otimes_{\bk} \bk[w_2^2, w_2^{-2}]$. Therefore, $\widehat{\Ho}^*(G)$ is of the form $\Lambda(x) \otimes_{\bk} \bk[z,z^{-1}]$, where $x$ and $z$ are of degrees $3$ and $4$, respectively, subject to the graded-commutative relations. 

By the decomposition Lemma~\ref{THdecomp}, $\widehat{\Ho}^*(G,\kG) \cong \widehat{\Ho}^*(G) \oplus \widehat{\Ho}^*(N)$ as graded $\bk$-modules. We then can define elements of the Tate-Hochschild cohomology ring of $\kG$ as follows. Since $\psi_1$ is an algebra monomorphism, we may identify any element of $\widehat{\Ho}^*(G)$ with its image under $\psi_1$. Let $E_i= \psi_i(1)$, $W_i=\psi_2(w_i)$, for $i=1,2$, and $W_2^{-1}=\psi_2(w_2^{-1})$. For simplification, we will use $C:= E_2+1$ in the following theorem. 

\begin{theorem}
Let $\bk$ be a field of characteristic $3$ and $S_3$ be the symmetric group on three elements. Then the Tate-Hochschild cohomology $\widehat{\HH}^*(\bk S_3,\bk S_3)$ of $S_3$ is generated as an algebra by elements $x, z, z^{-1}, C, W_1, W_2$, and $W_2^{-1}$ of degrees $3,4,-4, 0, 1, 2$, and $-2$, respectively, subject to the following relations:
\[xW_1=0, \;\; xW_2=zW_1, \;\; z^{-1}W_1=(xz^{-1})W_2^{-1}, \]
\[C^2=CW_2^{-1}=CW_i=0 \; (i=1,2), \]
\[W_2^2=zC, \;\; W_2^{-2}=z^{-1}C, \;\; W_1W_2=xC, \;\; W_1W_2^{-1}=xz^{-1}C, \]
together with the graded-commutative relations. In particular, the algebra monomorphism $\psi_1: \widehat{\Ho}^*(S_3) \rightarrow \widehat{\Ho}^*(S_3,\bk S_3)$ induces an isomorphism modulo radicals. 
\end{theorem}

\vspace{-1em}
\begin{proof}
$\widehat{\HH}^*(\kG,\kG) \cong \widehat{\Ho}^*(G,\kG)$ is a graded-commutative $\bk$-algebra whose underlying $\bk$-module is isomorphic to $\widehat{\Ho}^*(G) \oplus \widehat{\Ho}^*(N)$. Here, $\widehat{\Ho}^*(G)$ is a graded subalgebra of $\widehat{\HH}^*(\kG,\kG)$ generated by $x$, $z$ and $z^{-1}$. Additionally, $\psi_2(\widehat{\Ho}^*(N))$ is a graded $\widehat{\Ho}^*(G)$-submodule of $\widehat{\HH}^*(\kG,\kG)$ generated by $E_2$, $W_1$, $W_2$ and $W_2^{-1}$. This follows from the discussion after the proof of Theorem \ref{THdecomp product}. Moreover, we will check that these generators satisfy the following conditions: 
\begin{enumerate}
 \item action on $\psi_2(\widehat{\Ho}^*(N))$ as an $\widehat{\Ho}^*(G)$-module, and 
 \item every product in $\psi_2(\widehat{\Ho}^*(N))$ can be expressed as the sum of an element of $\widehat{\Ho}^*(G)$ and a $\widehat{\Ho}^*(G)$-linear combination of the images under $\psi_2$ of the generators of $\widehat{\Ho}^*(N)$.
\end{enumerate} 
Therefore, it is clear that $\widehat{\HH}^*(\kG,\kG)$ is generated as a $\bk$-algebra by $x, z, z^{-1}, E_2, W_1, W_2$, and $W_2^{-1}$, subject to these conditions. The first line of the relations in the statement of the theorem satisfies the first condition. The second and third lines satisfy the second condition. We will check each of them in detail. 

The restriction $\res^G_N: \widehat{\Ho}^*(G) \rightarrow \widehat{\Ho}^*(N)$, which sends $x \mapsto w_1w_2$, $z \mapsto w_2^2$, and $z^{-1} \mapsto w_2^{-2}$, is injective. We also observe that by graded-commutativity of the Tate cohomology ring, every element of odd degree has square $0$. In particular, $w_1w_1= -w_1w_1$ implies $w_1^2=0$. One can check that $\widehat{\Ho}^*(N)$ is an $\widehat{\Ho}^*(G)$-module with action via $\res^G_N$:
\begin{align*}
&x \cdot w_1 = w_1w_2w_1= w_1^2w_2=0, \\
&x \cdot w_2 = w_1w_2w_2= (-1)^2 w_2w_1w_2 = (-1)^2 w_2w_2w_1 = z \cdot w_1, \\
&x \cdot w_2^{-1} = w_1w_2w_2^{-1}= w_1, \\
&z \cdot w_2 = w_2^2w_2 = w_2^3, \\
&z \cdot w_2^{-1} = w_2^2 w_2^{-1}= w_2, \\
&z^{-1} \cdot w_1 = w_2^{-1}w_2^{-1}w_1= (-1)^{-2}w_2^{-1}w_1w_2^{-1} = (-1)^{-2}w_1w_2^{-1}w_2^{-1}= (xz^{-1}) \cdot w_2^{-1}, \\
&z^{-1} \cdot w_2 = w_2^{-1}w_2^{-1}w_2 = w_2^{-1}, \\
&z^{-1} \cdot w_2^{-1} = (w_2^{-1})^3.
\end{align*}
Therefore, as an $\widehat{\Ho}^*(G)$-module, $\widehat{\Ho}^*(N)$ is generated by $1, w_1, w_2$ and $w_2^{-1}$, subject to the relation $x \cdot w_1 = 0$, $x \cdot w_2 = z \cdot w_1$, and $z^{-1} \cdot w_1 = (xz^{-1}) \cdot w_2^{-1}$. By the isomorphism in Lemma~\ref{THdecomp} and mapping through $\psi_2$, we obtain the first line of the relations. 

To check the second and third lines of the relations, we recall the fact that the submodule of invariants $(\kG)^G$ is the center $Z(\kG)$ of the group algebra $\kG$, which is generated by conjugacy class representatives of $G$. Therefore, we may identify the degree-0 Tate-Hochschild cohomology with a quotient of $Z(\kG)$, as $\widehat{\HH}^0(\kG,\kG) \cong \widehat{\Ho}^0(G,\kG)$ is a quotient of $\Ho^0(G,\kG)$. Under this identification, $E_i$ corresponds to (a quotient of) the sum of the group elements conjugate to $g_i$. In particular, 
\[E_2^2= (a+a^{-1})^2= a^2 + 2 + a^{-2} = a^{-1} -1 + a = E_2 - 1 \]
in characteristic 3, which implies
\[C^2= (E_2 + 1)^2= E_2^2 + 2E_2 + 1 = 3E_2 = 0. \]

For the rest of the relations, we use the product formula in Theorem~\ref{THdecomp product}. Let $\alpha$ and $\beta$ be elements of $\widehat{\Ho}^*(N)$, we have:
\[\psi_2(\alpha) \smile \psi_2(\beta)= \psi_2(b^*(\alpha \beta)) + \psi_1(\cor^G_N(\alpha b^*(\beta))). \]
Recall that $b^*: \widehat{\Ho}^*(N) \rightarrow \widehat{\Ho}^*(^b N) = \widehat{\Ho}^*(N)$. By checking on the definition of $b^*$ and the degrees of $w_i$, we see that $b^*(w_2^{-1})=-w_2^{-1}$ and $b^*(w_i)=-w_i$, for $i=1,2$. Moreover, as there are no degree-$1, 2$ and $-2$ elements in $\widehat{\Ho}^*(G)$, we have $\cor^G_N(w_1)=\cor^G_N(w_2)=\cor^G_N(w_2^{-1})=0$. Similarly, by checking on the cochain level and using Lemma~\ref{THdecomp prop} (10), for all $n \in \mathbb{Z}$, we obtain: 
$$\cor^G_N(w_2^n)= \begin{cases}
                   0, &\quad n \text{ is odd} \\
                   -z^{n/2}, &\quad n \text{ is even.}
                   \end{cases} $$
Hence, using Lemma~\ref{THdecomp prop} (10) again, 
$$\cor^G_N(w_1w_2^n)= \begin{cases}
                   		-xz^{(n-1)/2}, &\quad n \text{ is odd} \\
                   		0, &\quad n \text{ is even.}
                   		\end{cases} $$                    

Let $\alpha=1$ and $\beta=w_1$, using the product formula in Theorem~\ref{THdecomp product}, we obtain:
\[\psi_2(1) \smile \psi_2(w_1)= \psi_2(b^*(w_1)) + \psi_1(\cor^G_N(b^*(w_1))) = \psi_2(-w_1)+0 = -W_1.\]
So $CW_1= (E_2+1)W_1 = E_2W_1 + W_1 = -W_1 + W_1 = 0$. Similarly, let $\alpha=1$ and $\beta=w_2$ or $w_2^{-1}$, we show that $CW_2=0=CW_2^{-1}$. This proves the second line of the relations.

Let $\alpha=\beta=w_2$, we have:
\begin{align*}
W_2^2= \psi_2(w_2) \smile \psi_2(w_2) &= \psi_2(b^*(w_2^2)) + \psi_1(\cor^G_N(w_2b^*(w_2))) \\
&= \psi_2(\res^G_N z \smile 1) + \psi_1(z) \\
&= z \smile \psi_2(1) + z \\
&= zE_2 + z = zC.
\end{align*}
Similarly, for $\alpha=\beta=w_2^{-1}$, we acquire that $W_2^{-2}= z^{-1}C$. 

Let $\alpha=w_1$ and $\beta=w_2^{-1}$:
\begin{align*}
W_1W_2^{-1}= \psi_2(w_1) \smile \psi_2(w_2^{-1}) &= \psi_2(b^*(w_1w_2^{-1})) + \psi_1(\cor^G_N(w_1b^*(w_2^{-1}))) \\
&= \psi_2(\res^G_N xz^{-1} \smile 1) + \psi_1(xz^{-1}) \\
&= xz^{-1} \smile \psi_2(1) + xz^{-1} \\
&= xz^{-1}E_2 + xz^{-1} = xz^{-1}C.
\end{align*}
Using the same argument, for $\alpha=w_1$ and $\beta=w_2$, we obtain $W_1W_2=xC$. Thus, we have found all necessary relations for the generators of $\widehat{\HH}^*(\kG,\kG)$. Furthermore, because the ring $\widehat{\Ho}^*(G,\kG)$ is graded-commutative, its nilpotent elements all lie in its radical. We observe that $C^2=0=W_1^2$, $W_2^3= W_2^2 W_2 = zCW_2=0$, and $(W_2^{-1})^3= W_2^{-2} W_2^{-1} = z^{-1}CW_2^{-1}=0$. This implies that $C, W_1, W_2$, and $W_2^{-1}$ are contained in the radical of the Tate-Hochschild cohomology ring $\widehat{\HH}^*(\kG,\kG)$. Consequently, modulo radicals, the algebra monomorphism $\psi_1: \widehat{\Ho}^*(S_3) \rightarrow \widehat{\HH}^*(\bk S_3,\bk S_3)$ induces an isomorphism. 
\end{proof}



\begin{thebibliography}{99}

\bibitem{AvMa} Luchezar L. Avramov and Alex Martsinkovsky, ``Absolute, relative, and Tate cohomology of modules of finite Gorenstein dimension,'' Proc. London Math. Soc. (3) \textbf{85} (2002), 393-440. 
\bibitem{BeCa} D. J. Benson and J. F. Carlson, ``Products in negative cohomology,'' J. Pure \& Applied Algebra, \textbf{82} (1992), 107-129.
\bibitem{BeJo} P. A. Bergh and D. A. Jorgensen, ``Tate-Hochschild homology and cohomology of Frobenius algebras,'' to appear in Journal of Noncommutative Geometry, (2011).
\bibitem{Ben} D. J. Benson, ``Commutative Algebra in the Cohomology of Groups,'' Trends in Commutative Algebra, MSRI Publications, \textbf{51} (2004), 1-50.
\bibitem{Ben1} D. J. Benson, \textit{Representations and Cohomology I: Basic representation theory of finite groups and associative algebras}, Cambridge Studies in Advanced Mathematics, Cambridge University Press, \textbf{30} (1991).
\bibitem{Ben2} D. J. Benson, \textit{Representations and Cohomology II: Cohomology of groups and modules}, Cambridge Studies in Advanced Mathematics, Cambridge University Press, \textbf{31} (1991).
\bibitem{Brown} K. S. Brown, \textit{Cohomology of Groups}, Graduate texts in mathematics, \textbf{87} Springer-Verlag, 1982.
\bibitem{CaEi} H. Cartan and S. Eilenberg, \textit{Homological Algebra}, Princeton University Press, Princeton, NJ, 1956.
\bibitem{CiSo} C. Cibils and A. Solotar, ``Hochschild cohomology algebra of abelian groups,'' Arch. Math, \textbf{68} (1997), 17-21.
\bibitem{Nguyen} V. C. Nguyen, ``Tate and Tate-Hochschild Cohomology for finite dimensional Hopf Algebras,'' (submitted) arXiv:1209.4888.
\bibitem{SW} S. F. Siegel and S. J. Witherspoon, ``The Hochschild cohomology ring of a group algebra,'' Proc. London Math. Soc. \textbf{79} (1999), 131-157. 
\bibitem{Tate} J. Tate, ``The higher dimensional cohomology groups of class field theory,''  Ann. of Math. (2) \textbf{56} (1952), 294-297.


\end{thebibliography}
\end{document}